\documentclass[graybox]{article}

\usepackage{amsmath}
\usepackage{hyperref}
\numberwithin{equation}{section}
\usepackage{cleveref}
\usepackage{mathptmx}       
\usepackage{helvet}         
\usepackage{courier}        
\usepackage{type1cm}        
\usepackage{amssymb}
\usepackage{makeidx}         
\usepackage{graphicx}        
\usepackage{multicol}        
\usepackage[bottom]{footmisc}

\makeindex             

\begin{document}
	
	\title{Modeling and Simulation of Macroscopic Pedestrian Flow Models}
	\author{Naveen Kumar Mahato  and Axel Klar and Sudarshan Tiwari\\ Technische Universität Kaiserslautern,\\ Erwin-Schrödinger-Straße, 67663 Kaiserslautern, Germany\\ mahato@mathematik.uni-kl.de,\\ klar@mathematik.uni-kl.de, \\tiwari@mathematik.uni-kl.de}
	%
	%
	\maketitle

\abstract{We analyze numerically some macroscopic models of pedestrian motion such as Hughes model (Hughes in Transportation research Part B: Methodological 36: 507-535, 2002) and mean field game with nonlinear mobilities (Burger et al. in Discrete and Continuous Dynamical Systems. Series B. A Journal Bridging Mathematics and Sciences 19(5): 1311--1333, 2014) modeling fast exit scenarios in pedestrian crowds. A model introduced by Hughes consisting of a non-linear conservation law for the density of pedestrians coupled with an Eikonal equation for a potential modeling the common sense of the task. Mean field game with nonlinear mobilities is obtained by an optimal control approach, where the motion of every pedestrian is determined by minimizing a cost functional, which depends on the position, velocity, exit time and the overall density of people. We consider a parabolic optimal control problem of nonlinear mobility in pedestrian dynamics, which leads to a mean field game structure. We show how optimal control problem related to the Hughes model for pedestrian motion.  Furthermore we provide several numerical results which relate both models in one and two dimensions.}

\section{Introduction}
Mathematical modeling and numerical simulation of human crowd motion have become a major subject of research with a wide field of applications. A variety of models for pedestrian behavior have been proposed on different levels of description in recent years. Macroscopic pedestrian flow model involving equations for density and mean velocity of the flow are derived in Refs. ~\cite{bellomo, burger, hughes, mahato, mahato2}.

\section{Optimal control problem of pedestrian flow from \protect\cite{burger}}
For completeness of the presentation up to higher dimensions we review the macroscopic optimal control problem for pedestrian flow, see Refs. ~\cite{burger, burger2, burger3}. For further details we refer to these papers. There, denoting the (normalized) density function of the pedestrians by $\rho(t,x)$ and the momentum (or the flux density) by $m = F(\rho) v$ at position $x\in\Omega$, velocity $v\in\Omega$ and time $t$, where the function $F(\rho(t,x))$ describing the nonlinear mobility of the pedestrians (or the costs created by large densities) and $\Omega\in\mathbb{R}^d$, $d = 1,2$ is a bounded domain representing the pedestrian area. We assume the boundary $\partial\Omega$ is split into a Neumann part $\Gamma_N \subseteq \partial\Omega$ modeling walls or obstacles, $\Gamma_E \subseteq \partial\Omega$ modeling the exits such that $\partial\Omega = \Gamma_N\bigcup\Gamma_E$ and $\Gamma_N\bigcap\Gamma_E = \phi$. If we denote the rate of passing the exit by $\beta$, then we have an outflow proportional to $\beta\rho$. Hence, for a stochastic particles and a final time $T$ sufficiently large, the minimization functional is given by the following parabolic optimal control problem: 
\begin{subequations}\label{eq1}
	\begin{equation}\label{eq1a}
		\underset{(\rho,m)}{\min} \hspace{1mm} I_T(\rho,m) = \underset{(\rho,m)}{\min} \hspace{2mm} \frac{1}{2}\int_{0}^{T}\int_{\Omega} \dfrac{\lvert m(t,x)\rvert^2}{F(\rho(t,x))}  dxdt + \frac{\alpha}{2} \int_{0}^{T}\int_{\Omega} \rho(t,x)dxdt,
	\end{equation}
	subject to
	\begin{align}
		\partial_t \rho + \nabla \cdot m &= \frac{\sigma^2}{2} \Delta\rho,\hspace{10mm} \mbox{in} \hspace{2mm} \Omega\times(0,T), \label{eq1b}\\
		\left( m - \frac{\sigma^2}{2}\nabla\rho \right) \cdot n &= 0,\hspace{17mm} \mbox{on} \hspace{2mm} \Gamma_N\times(0,T), \label{eq1c}\\
		\left( m - \frac{\sigma^2}{2}\nabla\rho \right) \cdot n &= \beta\rho,\hspace{14.5mm} \mbox{on} \hspace{2mm} \Gamma_E\times(0,T), \label{eq1d}\\
		\rho(0,x) &= \rho_0(x), \hspace{12mm} \mbox{in} \hspace{2mm} \Omega.\label{eq1e}
	\end{align}  
\end{subequations} 
This optimality system can be seen as the mean field games structure, see Ref. \cite{LL}. We start with defining the Lagrangian with dual variable $\Phi = \Phi(t,x)$ as
\begin{align*}
	 &L_T(\rho,m,\Phi) = I_T(\rho,m)+\int_{0}^{T}\int_{\Omega}(\partial_t \rho + \nabla \cdot m -  \frac{\sigma^2}{2} \Delta\rho)\Phi dx dt \\
	& = I_T(\rho,m) + \int_{0}^{T}\int_{\Omega} \left[ \rho \left( -\partial_t \Phi -  \frac{\sigma^2}{2} \Delta\Phi\right) - m \cdot \nabla\Phi\right]   dx dt \\
	& + \int_{0}^{T}\int_{\partial\Gamma_E} \left( \underbrace{-\frac{\sigma^2}{2} \nabla\rho \cdot n\Phi + m \cdot n\Phi }_{\beta\rho\Phi} + \frac{\sigma^2}{2} \rho \nabla\Phi \cdot n  \right)  ds dt  + \int_{0}^{T}\int_{\partial\Gamma_N}  \frac{\sigma^2}{2} \rho \nabla\Phi \cdot n  ds dt .
\end{align*}
The optimality condition with respect to $m$ and $\rho$, yields the following equations  
\begin{align*}
0 &=  \partial_m L_T(\rho,m,\Phi) = \dfrac{m(t,x)}{F(\rho(t,x))} - \nabla\Phi, \\
 0 &=  \partial_\rho L_T(\rho,m,\Phi) =	- \frac{1}{2} \dfrac{\lvert m(t,x)\rvert^2F'(\rho)}{F^2(\rho)}+\frac{\alpha}{2}-\partial_t \Phi - \frac{\sigma^2}{2} \Delta\Phi.
\end{align*} 
Inserting $ m = F(\rho(t,x))\nabla\Phi$ we obtain the following system of equations
\begin{subequations}\label{eq2}
	\begin{align}
		\partial_t \rho + \nabla \cdot (F(\rho) \nabla\Phi) -  \frac{\sigma^2}{2} \Delta\rho &= 0,\hspace{31.5mm} \mbox{in} \hspace{2mm} \Omega\times(0,T), \label{eq2a}\\
		\partial_t \Phi + \frac{F'(\rho)}{2}\lvert\nabla\Phi\rvert^2 +  \frac{\sigma^2}{2} \Delta\Phi &= \frac{\alpha}{2}, \hspace{30mm} \mbox{in} \hspace{2mm} \Omega\times(0,T), \label{eq2b}\\
		\left(F(\rho) \nabla\Phi - \frac{\sigma^2}{2}\nabla\rho \right) \cdot n = 0,  & \hspace{3mm} \frac{\sigma^2}{2}\nabla\Phi \cdot n = 0, \hspace{15mm} \mbox{on} \hspace{2mm} \Gamma_N\times(0,T), \label{eq2c}\\
		\left(F(\rho) \nabla\Phi - \frac{\sigma^2}{2}\nabla\rho \right) \cdot n = \beta\rho, &\hspace{3mm} \frac{\sigma^2}{2}\nabla\Phi \cdot n + \beta\rho = 0, \hspace{6.5mm} \mbox{on} \hspace{2mm} \Gamma_E\times(0,T), \label{eq2d}\\
		\rho(0,x) = \rho_0(x), & \hspace{3mm} \Phi(T,x) = 0, \hspace{19mm} \mbox{in} \hspace{2mm} \Omega.\label{eq2e}
	\end{align}  
\end{subequations}
System (\ref{eq2}) has the structure of a mean field games for pedestrian dynamics, which contains the Fokker-Planck equation (\ref{eq2a}) has to be solved forward in time and the Hamilton-Jacobi equation (\ref{eq2b}) that has to be solved backward in time.
\section{Relation to the classical Hughes model \protect\cite{hughes}}
In this Section we discuss the relation which shows that for vanishing viscosity $\sigma = 0$ of the optimality system (\ref{eq1}) has a similar structure as the classical Hughes model for pedestrian flow. Hughes proposed that pedestrians seek the fastest path to the exit, but at the same time try to avoid congested areas, for details see Ref. ~\cite{hughes}. Let us consider the governing equations of Hughes model for pedestrian flow,
\begin{subequations}\label{eq3}
	\begin{align}
		\partial_t \rho - \nabla\cdot(\rho f^2(\rho)\nabla\Phi) -  \frac{\sigma^2}{2} \Delta\rho &= 0, \hspace{25mm} \mbox{in} \hspace{2mm} \Omega\times(0,T), \label{eq3a}\\
		\lvert\nabla\Phi\lvert &= \frac{1}{f(\rho)}, \hspace{19mm} \mbox{in} \hspace{2mm} \Omega\times(0,T), \label{eq3b}\\
		\left(\rho f^2(\rho) \nabla\Phi - \frac{\sigma^2}{2}\nabla\rho \right) \cdot n = 0,  & \hspace{3mm} \Phi = \infty, \hspace{18.5mm} \mbox{on} \hspace{2mm} \Gamma_N\times(0,T), \label{eq3c}\\
		\left(\rho f^2(\rho) \nabla\Phi - \frac{\sigma^2}{2}\nabla\rho \right) \cdot n = \beta\rho, &\hspace{3mm} \Phi = 0, \hspace{19.5mm} \mbox{on} \hspace{2mm} \Gamma_E\times(0,T), \label{eq3d}\\
		\rho(0,x) &= \rho_0(x), \hspace{20mm} \mbox{in} \hspace{2mm} \Omega,\label{eq3e}
	\end{align}
\end{subequations}    
where the function $f(\rho) = \rho_{max} - \rho$ with $\rho_{max}$ denote the maximum density and models how pedestrians change their direction and velocity due to the surrounding density, i.e. provides a weighting or cost with respect to high densities. Saturation effects are included via the function $f(\rho)$ for $\rho \longrightarrow \rho_{\max}$. 

On the other hand, if we choose the mobility/penalization function for high densities such as $F(\rho)=\rho f(\rho)^2$, then the optimality system (\ref{eq2}) for vanishing viscosity can be written as
\begin{subequations}\label{eq4}
	\begin{align}
		\partial_t \rho + \nabla \cdot (\rho f(\rho)^2 \nabla\Phi) &= 0,\hspace{11.5mm} \mbox{in} \hspace{2mm} \Omega\times(0,T), \label{eq4a}\\
		\partial_t \Phi + \frac{f(\rho)}{2}(f(\rho)+2\rho f'(\rho))\lvert\nabla\Phi\rvert^2 &= \frac{\alpha}{2}, \hspace{10mm} \mbox{in} \hspace{2mm} \Omega\times(0,T), \label{eq4b}
	\end{align}  
\end{subequations}
where the initial, terminal and boundary conditions are same as in system (\ref{eq2}). Furthermore, one can expect the equilibration of $\Phi$ backward in time for large $T$. Then for time $t$ of order one the limiting model becomes
\begin{subequations}\label{eq5}
	\begin{align}
		\partial_t \rho + \nabla \cdot (\rho f(\rho)^2 \nabla\Phi) &= 0,\hspace{31mm} \mbox{in} \hspace{2mm} \Omega\times(0,T), \label{eq5a}\\
		(f(\rho)+2\rho f'(\rho))\lvert\nabla\Phi\rvert^2 &= \frac{\alpha}{f(\rho)}, \hspace{25mm} \mbox{in} \hspace{2mm} \Omega\times(0,T). \label{eq5b}
	\end{align}          
\end{subequations}          
Hence, if we set $\alpha = 1$, the system (\ref{eq5}) is almost equivalent to the Hughes model (\ref{eq3}) for vanishing viscosity. Note that the sign difference in equations (\ref{eq3a}) and (\ref{eq5a}) is not an actual, since due to the signs in the backward equation we shall obtain $\Phi$ as the negative of the distance function used in the Hughes model.


\section{Numerical results}
In this Section we present a series of numerical experiments for the equations from both proposed models. We compare the relation between the models for different parameters in one and two dimensions. We use finite difference scheme for solving the classical Hughes model, where central difference in space and the forward difference in time, i.e. forward time centered space (FTCS) scheme for the nonlinear conservation law and an upwind Godunov scheme for the Eikonal equation.  We follow the steepest descent algorithm from Ref. ~\cite{burger}, to solve the mean field games structure, in which we use FTCS finite difference scheme to solve both forward and backward equations.

We consider a one-dimensional domain $\Omega = [-1, 1]$ with exits located at $x = \pm 1$ for the numerical simulation as a configuration defined in Ref. ~\cite{burger}. We choose the maximum density $\rho_{max}$, the weighting parameter $\alpha$ and the flow rate parameter $\beta$ as $1$. Furthermore, we consider the time interval as $t\in[0, 3]$. The time step is set to  $\Delta t = 10^{-4}$ for Hughes and $\Delta t = 10^{-3}$ for MFG. We use the spatial discretization $h = 10^{-2}$, the diffusion coefficient $\frac{\sigma^2}{2} = h$ and the initial density $\rho_0 = \frac{1}{3}$ in both models. 

\begin{figure}
	\includegraphics[width=.5\linewidth]{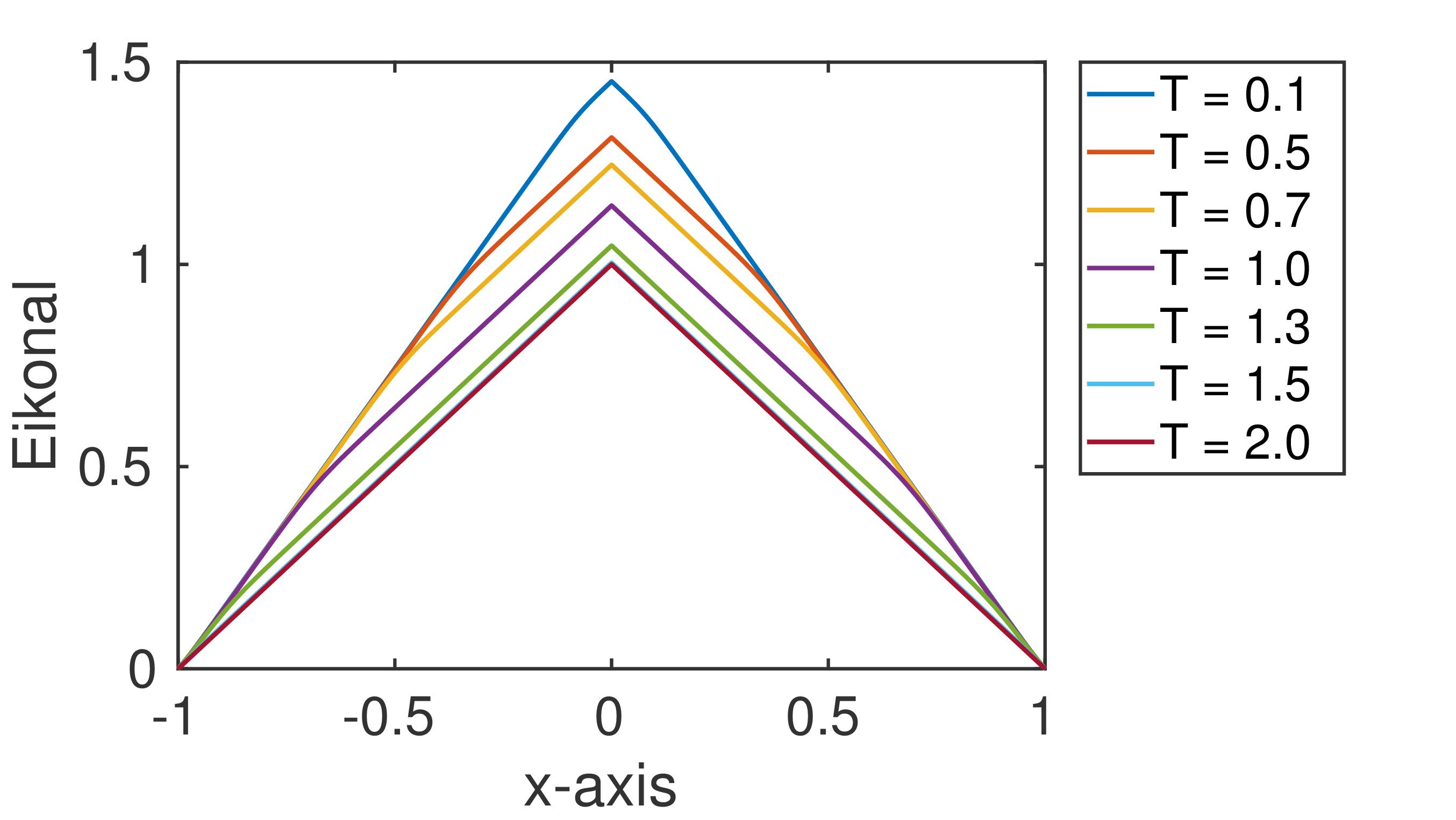}\hfill
	\includegraphics[width=.5\linewidth]{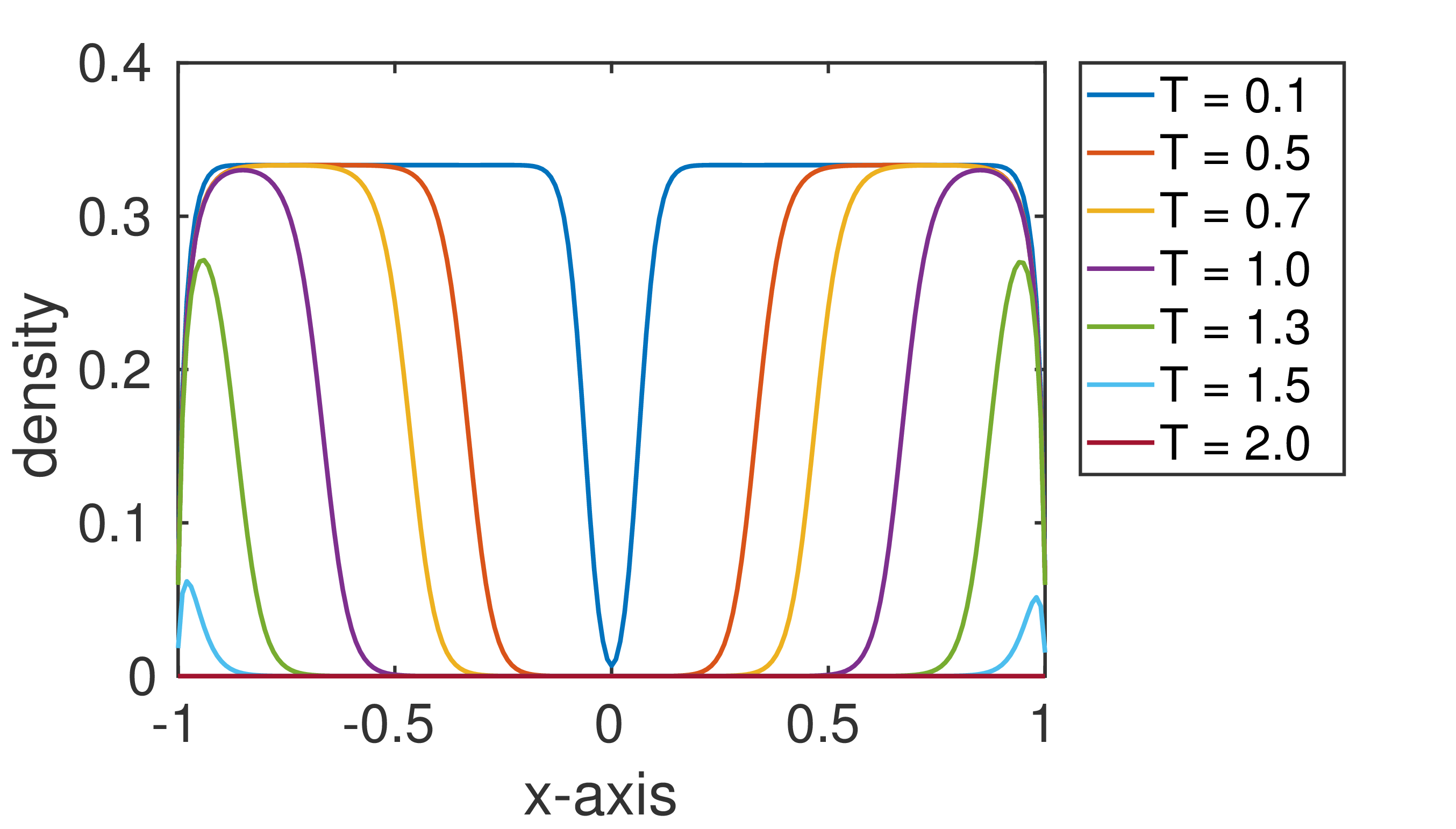}\\ 
	\includegraphics[width=.5\linewidth]{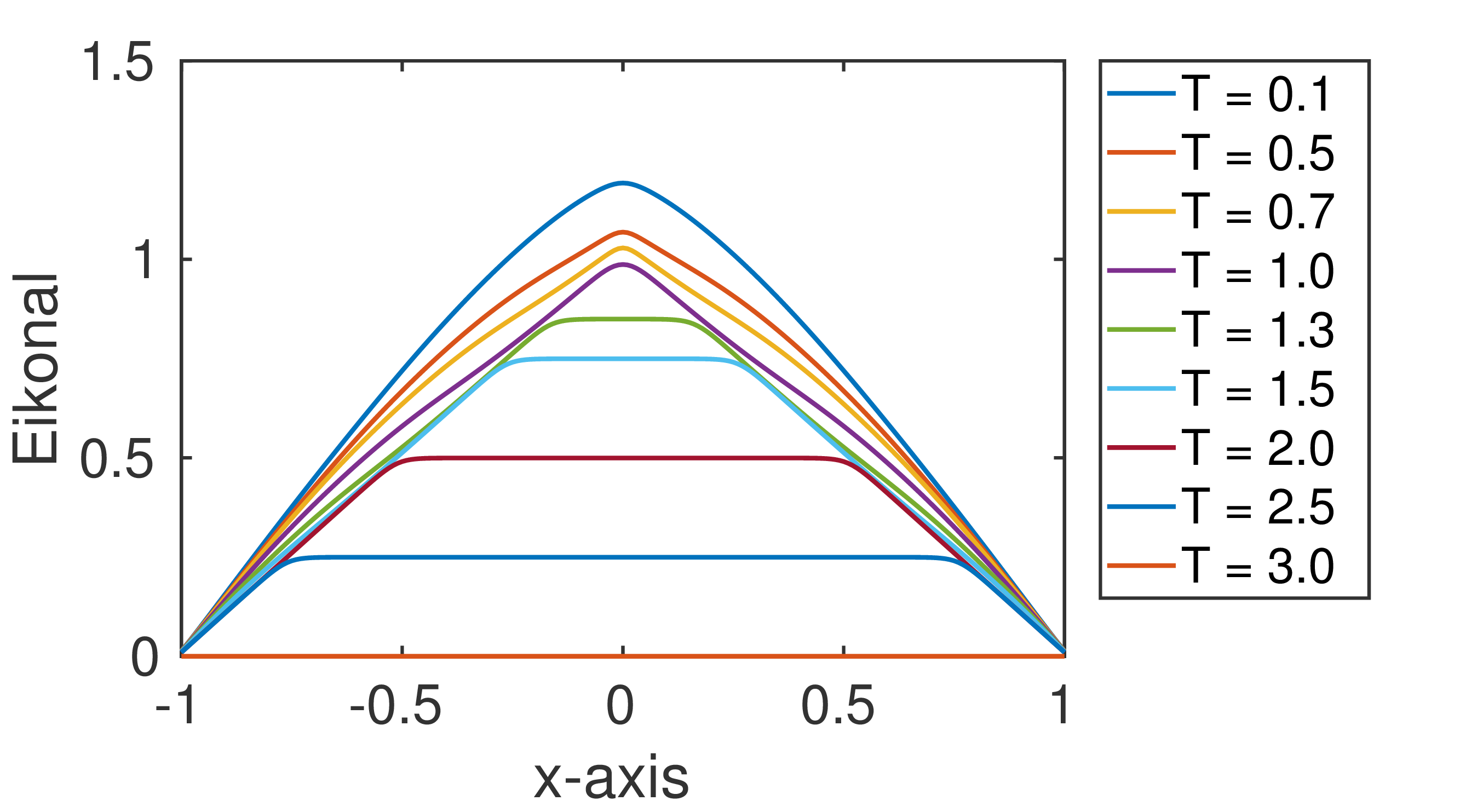}\hfill
	\includegraphics[width=.5\linewidth]{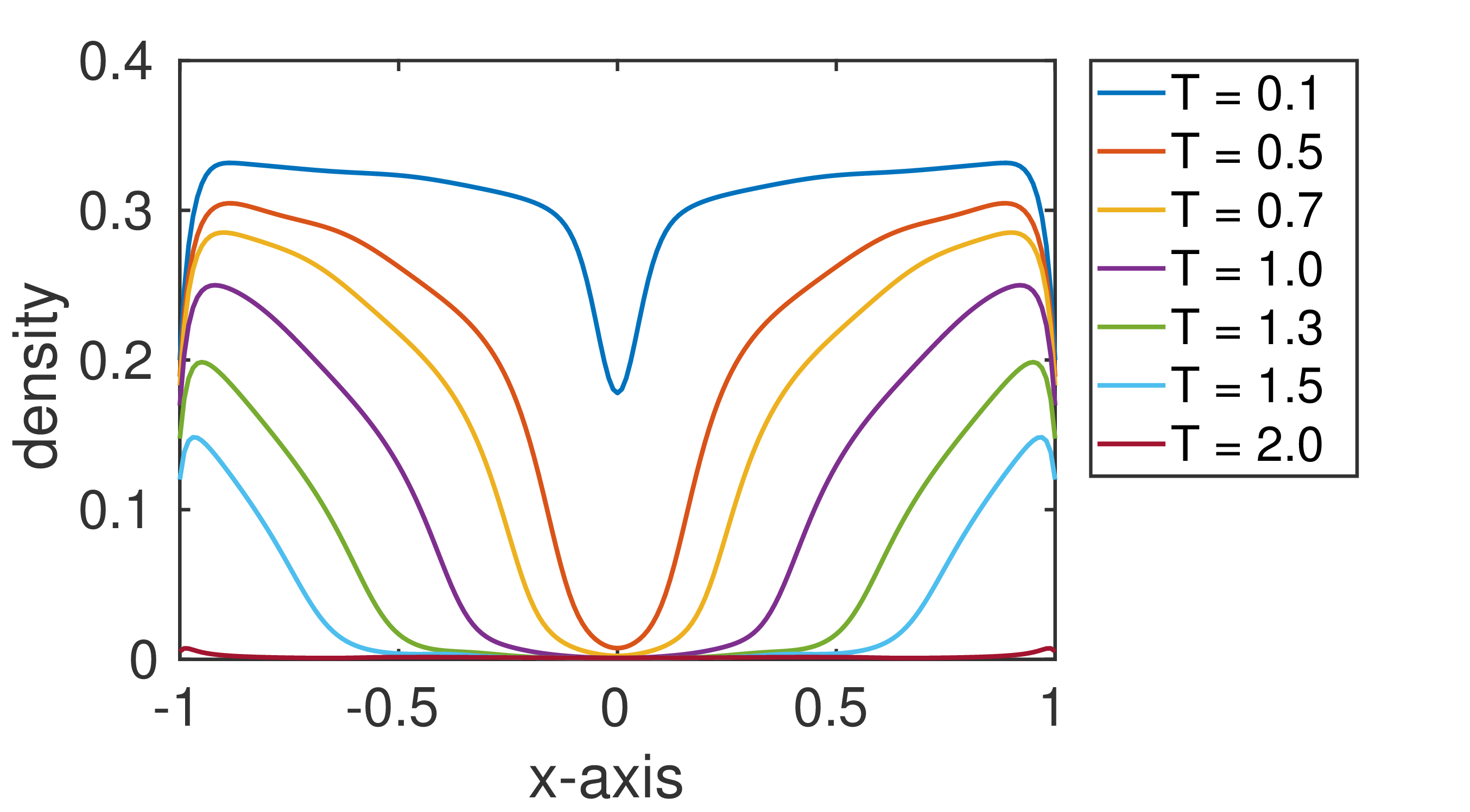}
	\caption{Evolution of solutions at different times for the Hughes model (top) and for the MFG structure (bottom).}\label{fig1}	
\end{figure}        
Figure \ref{fig1} shows the evolution of solutions at different times for both mean field type structure. One observes that the non-stationary Eikonal solution of the MFG structure has a similar behavior as the stationary Eikonal solution of the classical Hughes model until the density is not zero, as we expected equilibration of $\Phi$ in the equation (\ref{eq4b}). One also observes from the density solution that both models have similar behavior as pedestrians start in immediate vacuum formation at the center $x = 0$. Although the models have a very similar structure, pedestrians wait for a little while at the center and then start to move at a higher speed in the case of the mean field games compare to the Hughes model.

The extension of the above method into higher dimensions is straight forward. Here, we restrict ourselves to two-dimensional problems. Suppose the geometry for numerical experiment is taken as $\Omega = [-1, 1]\times [-1, 1]$ with exits located at $(\pm 1, \pm 1)$. Furthermore, we choose all parameters as for one dimension.
\begin{figure}
	\includegraphics[scale=0.03]{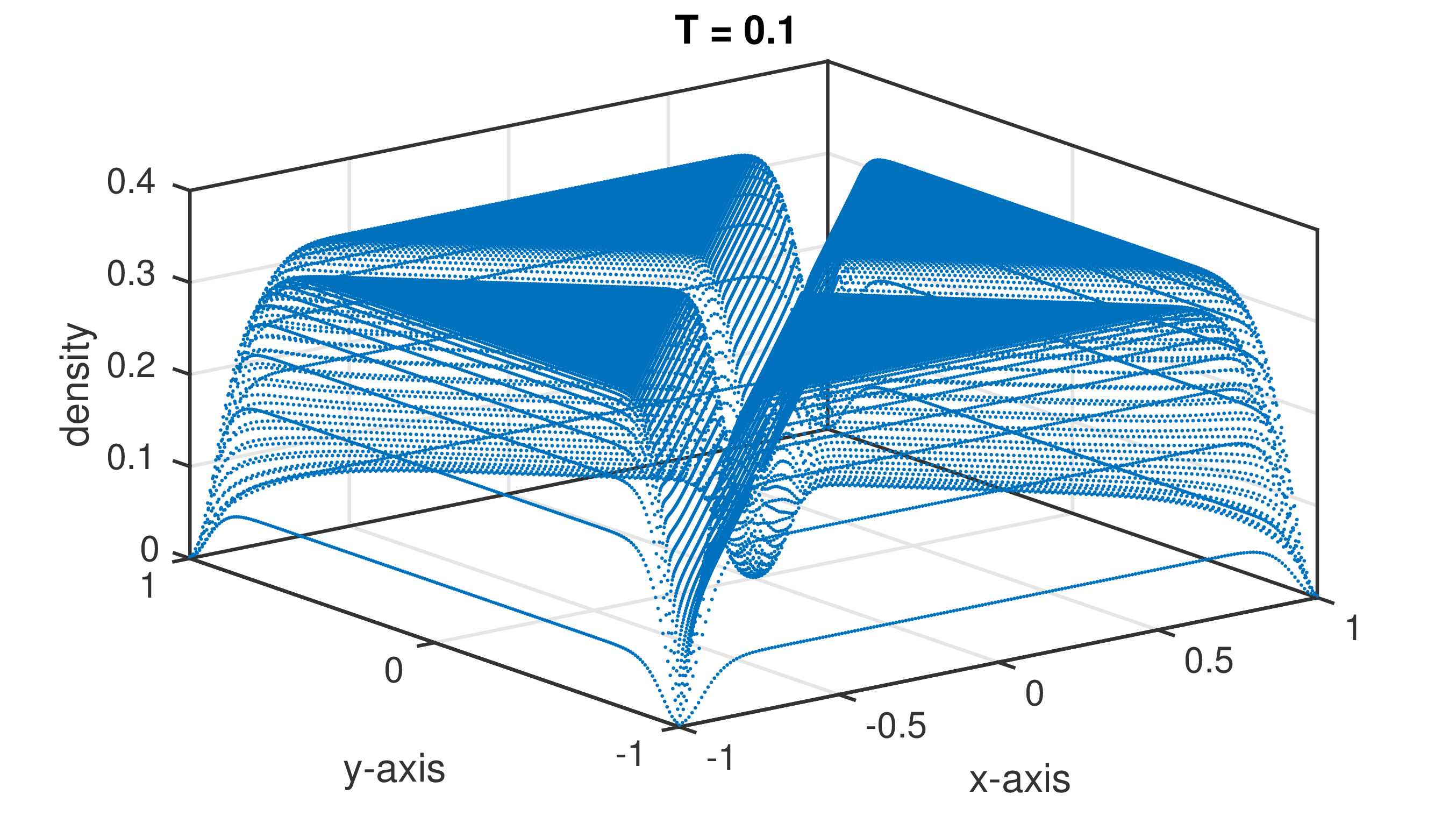}\hfill
	\includegraphics[scale=0.03]{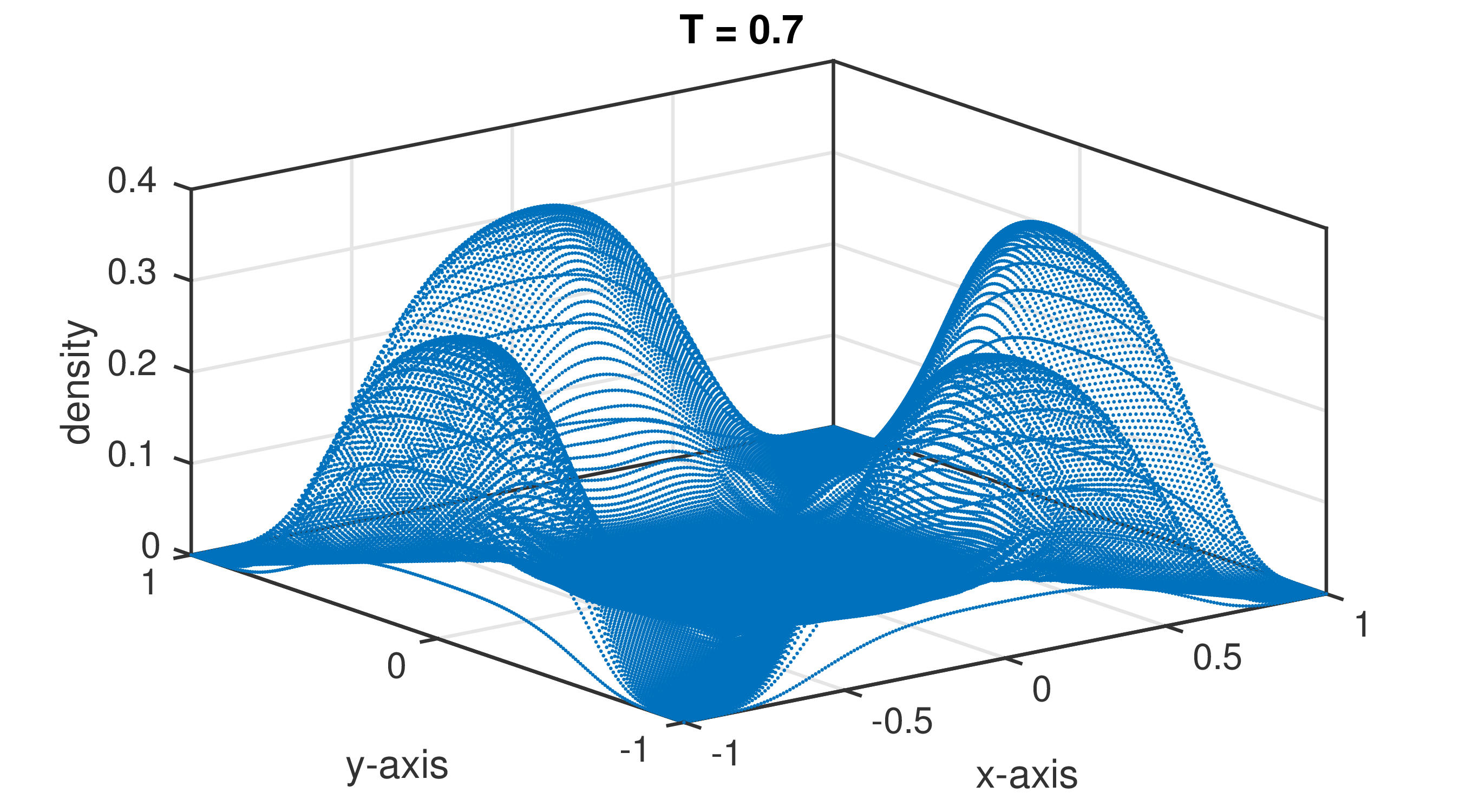}\hfill
	\includegraphics[scale=0.03]{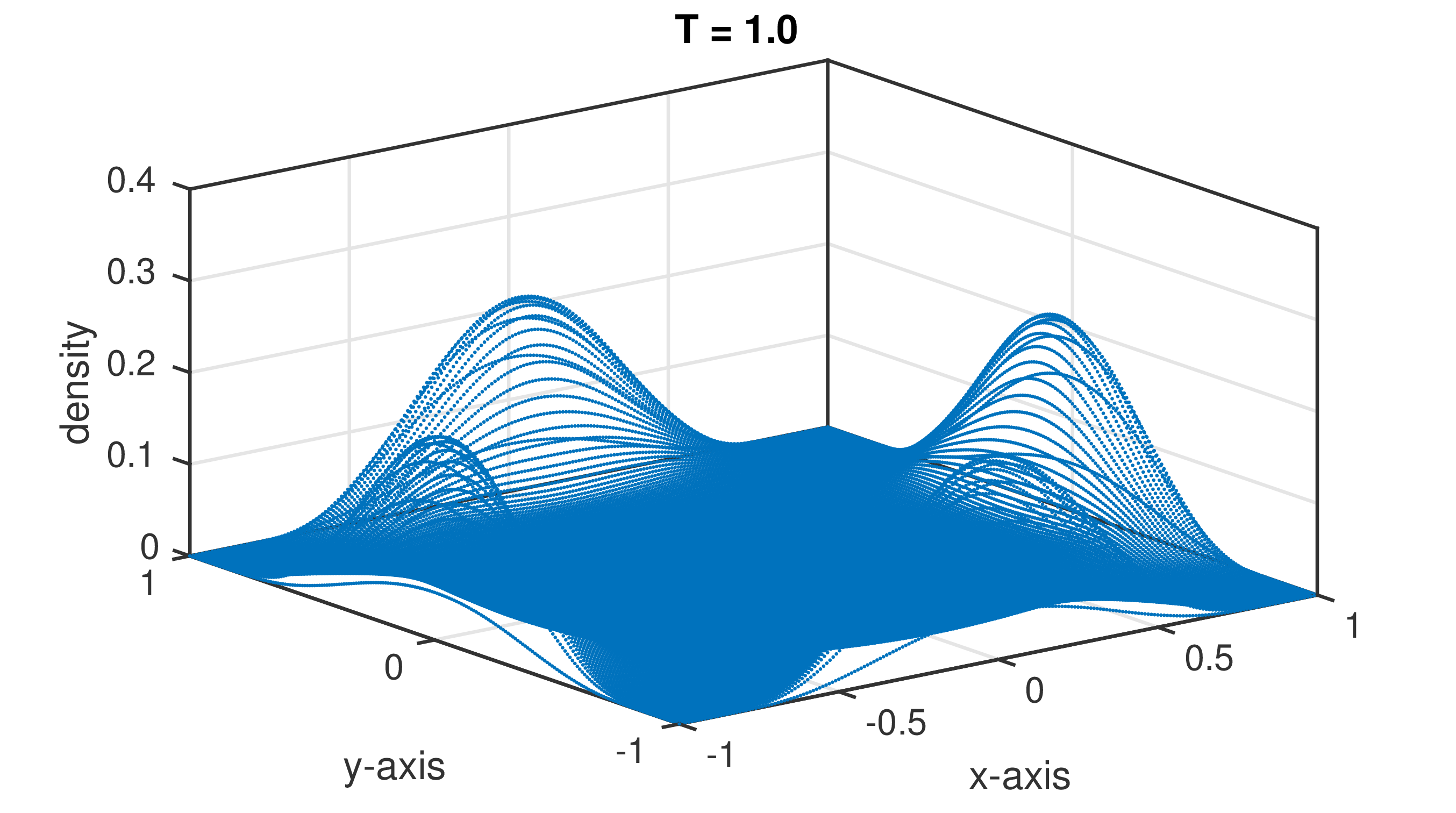}\hfill
	\includegraphics[scale=0.03]{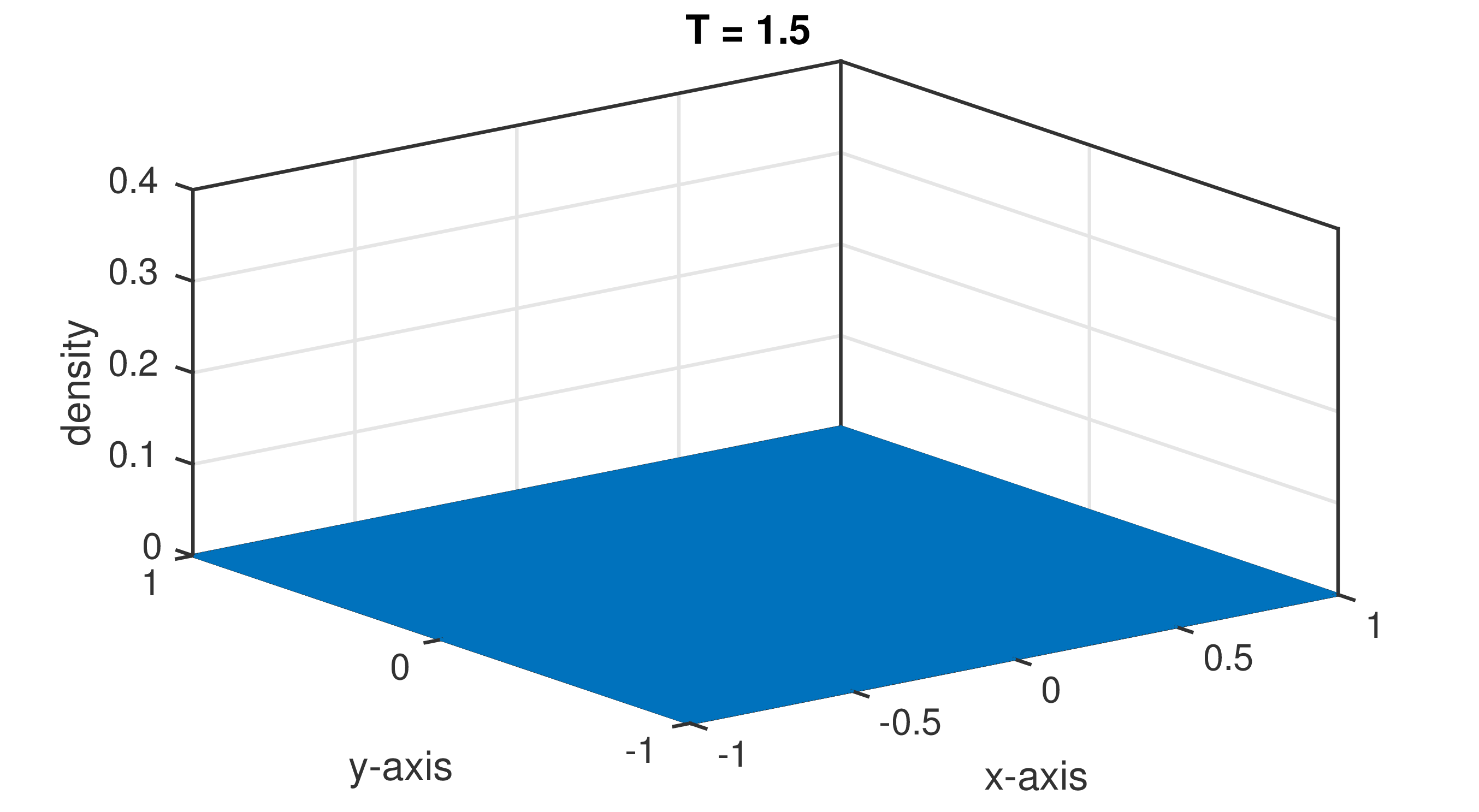}\\
	\includegraphics[scale=0.03]{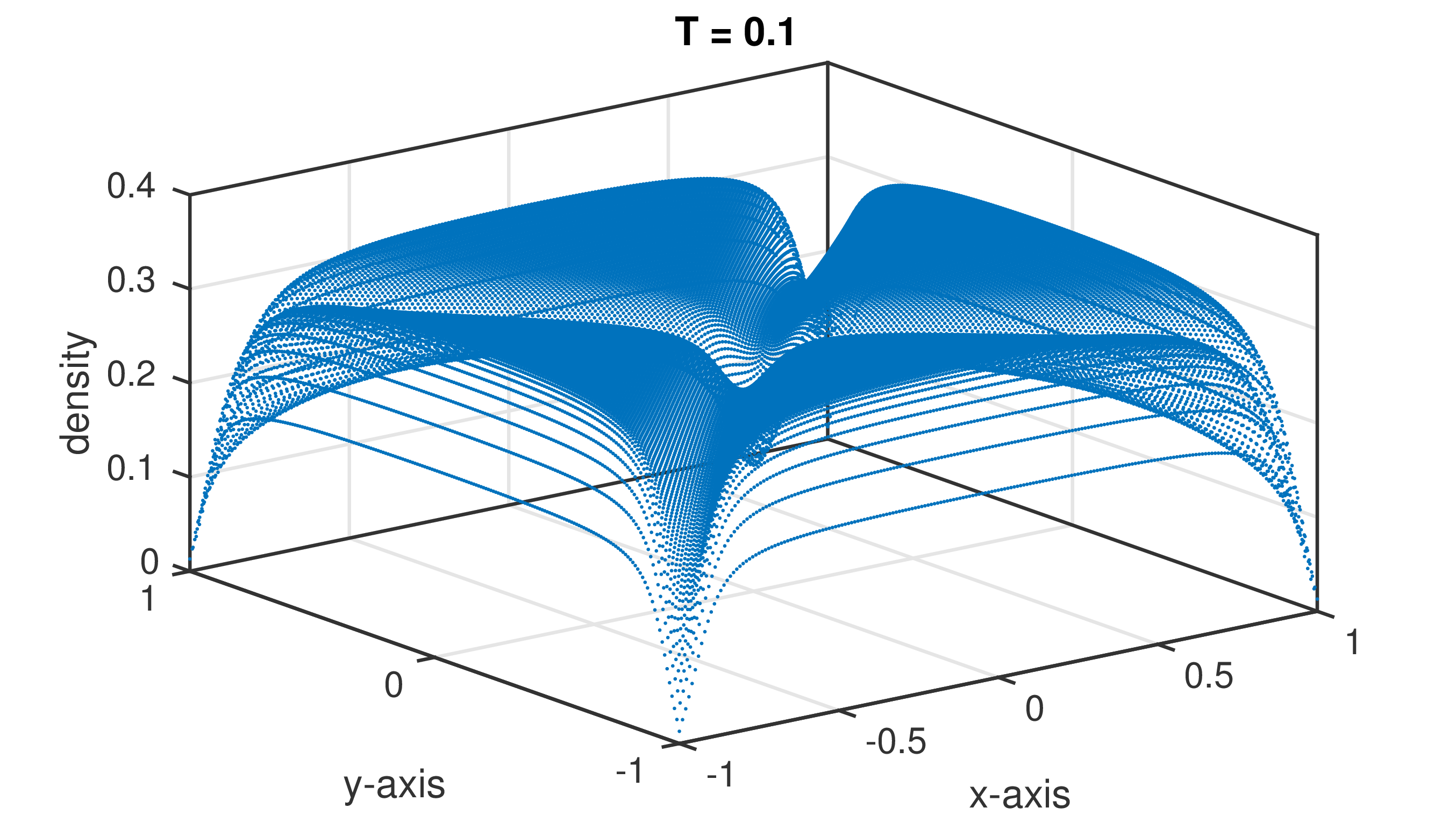}\hfill
	\includegraphics[scale=0.03]{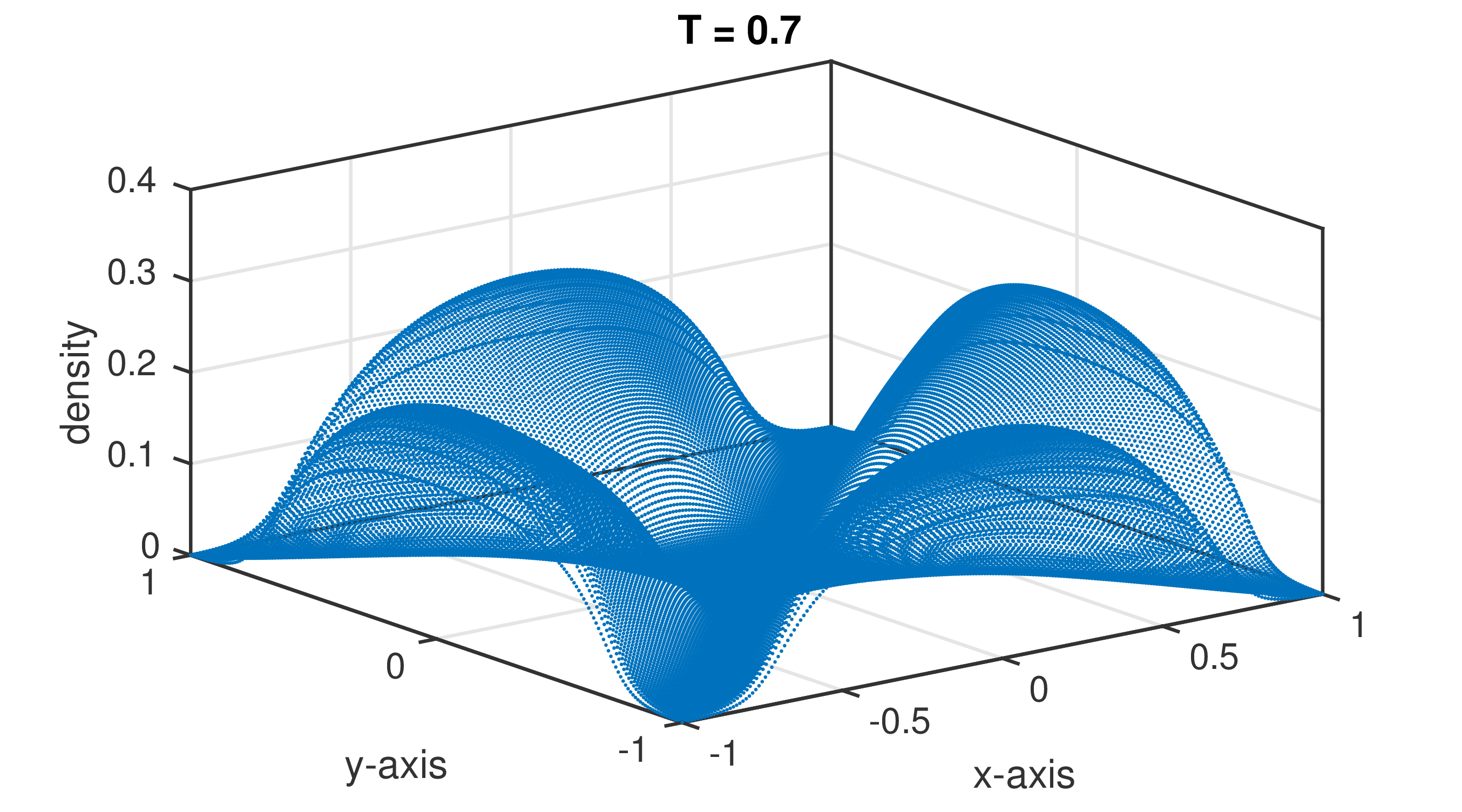}\hfill
	\includegraphics[scale=0.03]{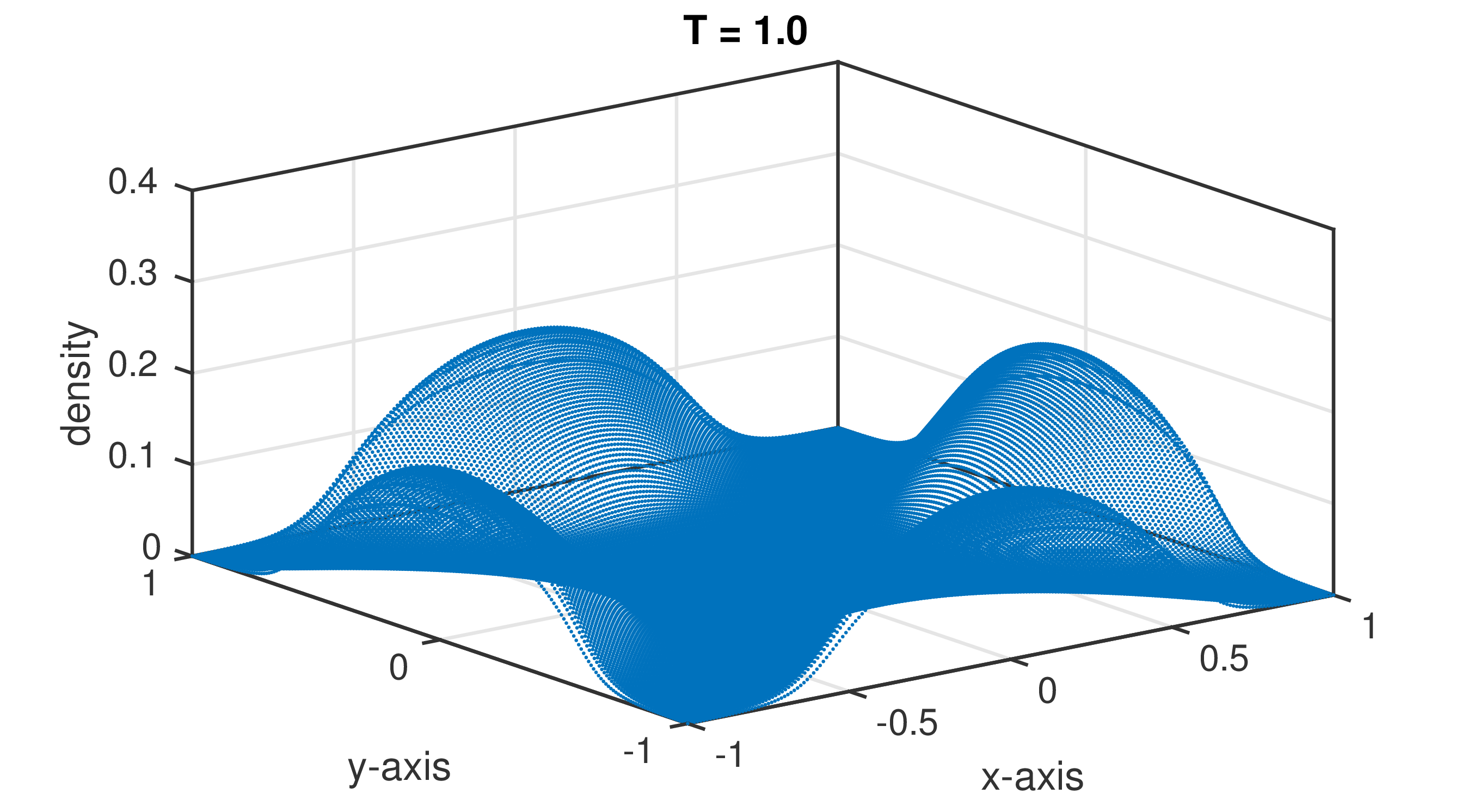}\hfill
	\includegraphics[scale=0.03]{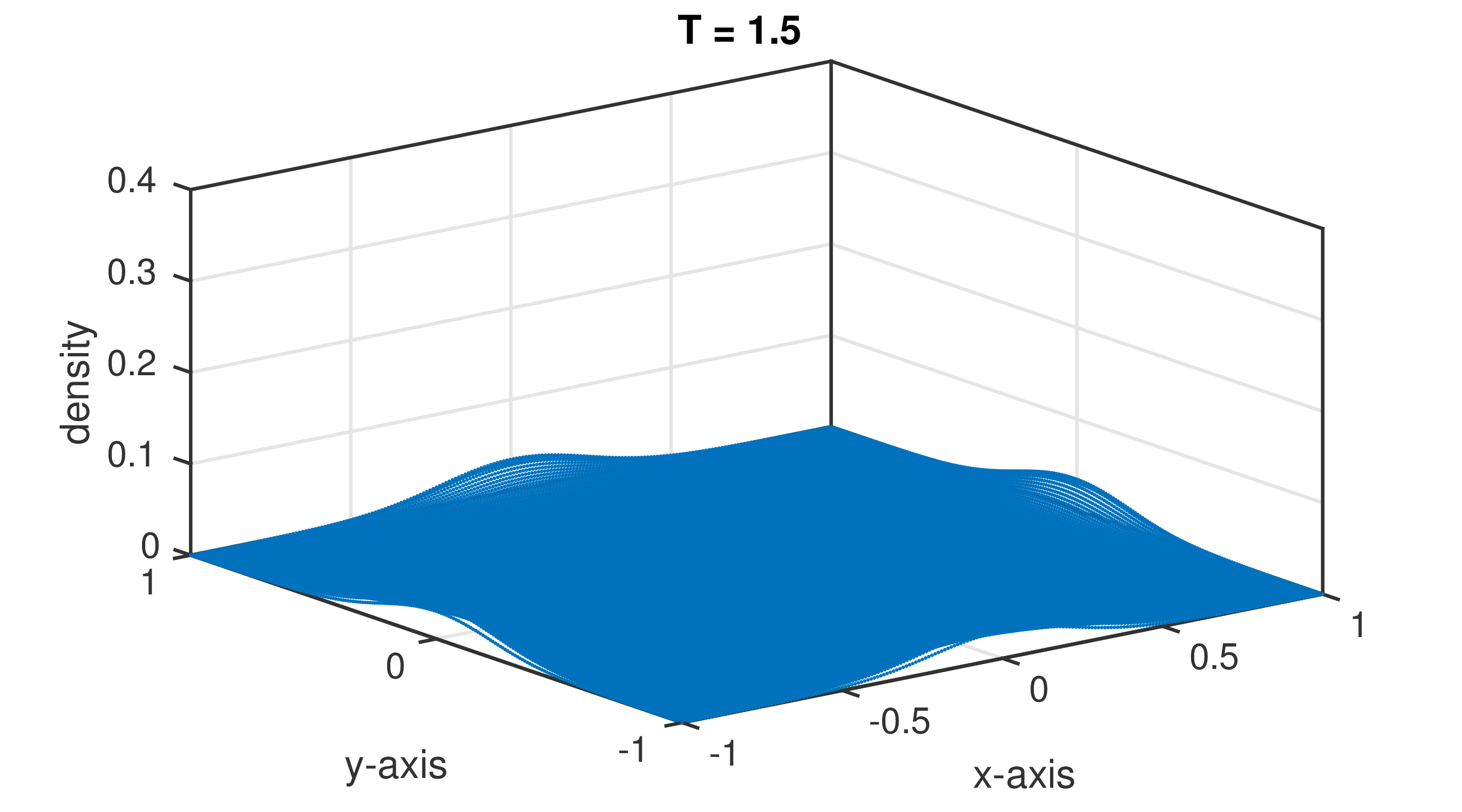}
	\caption{Evolution of the density solution at different times for the Hughes model (top) and the MFG structure (bottom).}\label{fig2}
	
\end{figure}     
\begin{figure}
	\includegraphics[width=.5\linewidth]{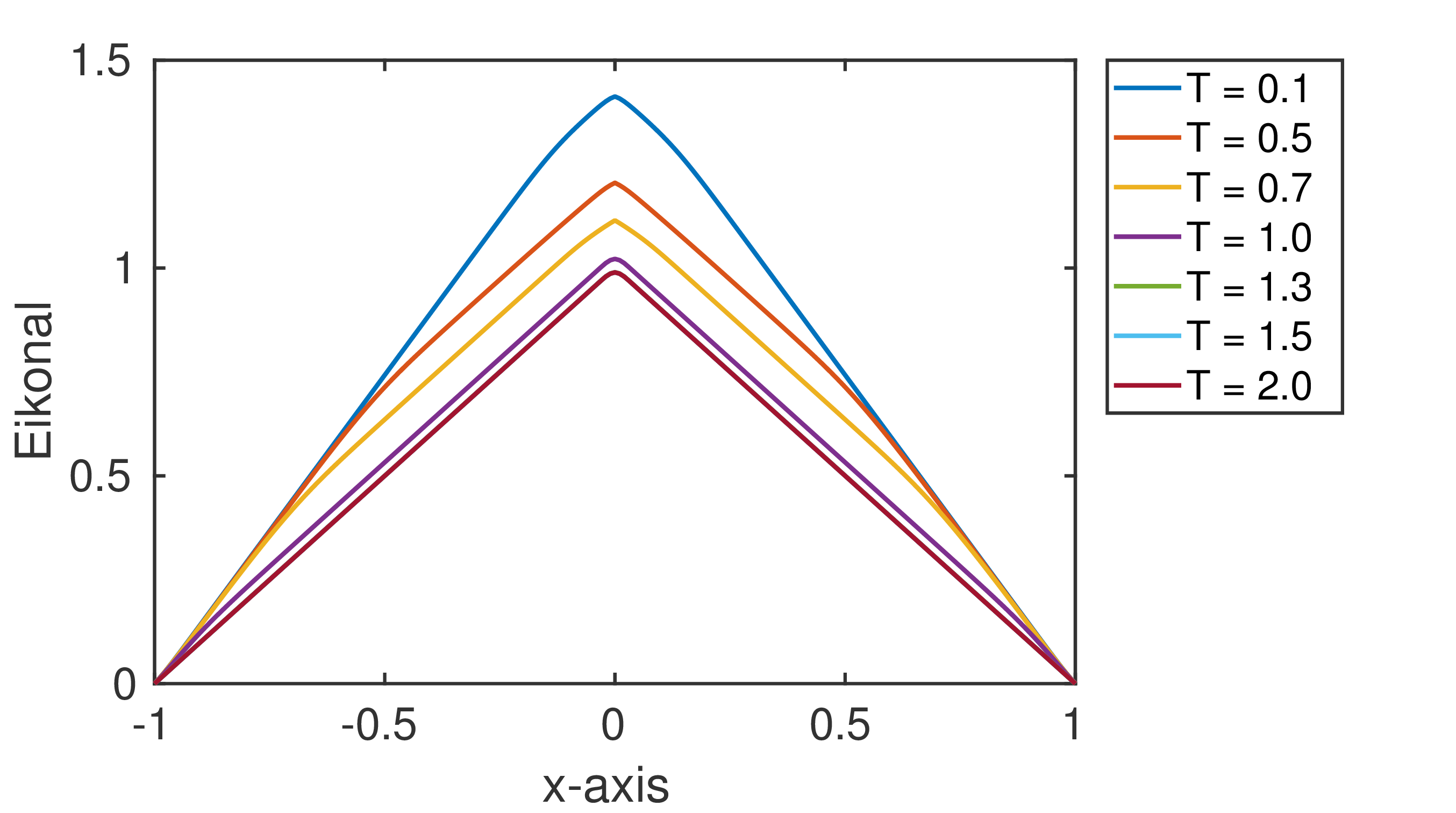}\hfill
	\includegraphics[width=.5\linewidth]{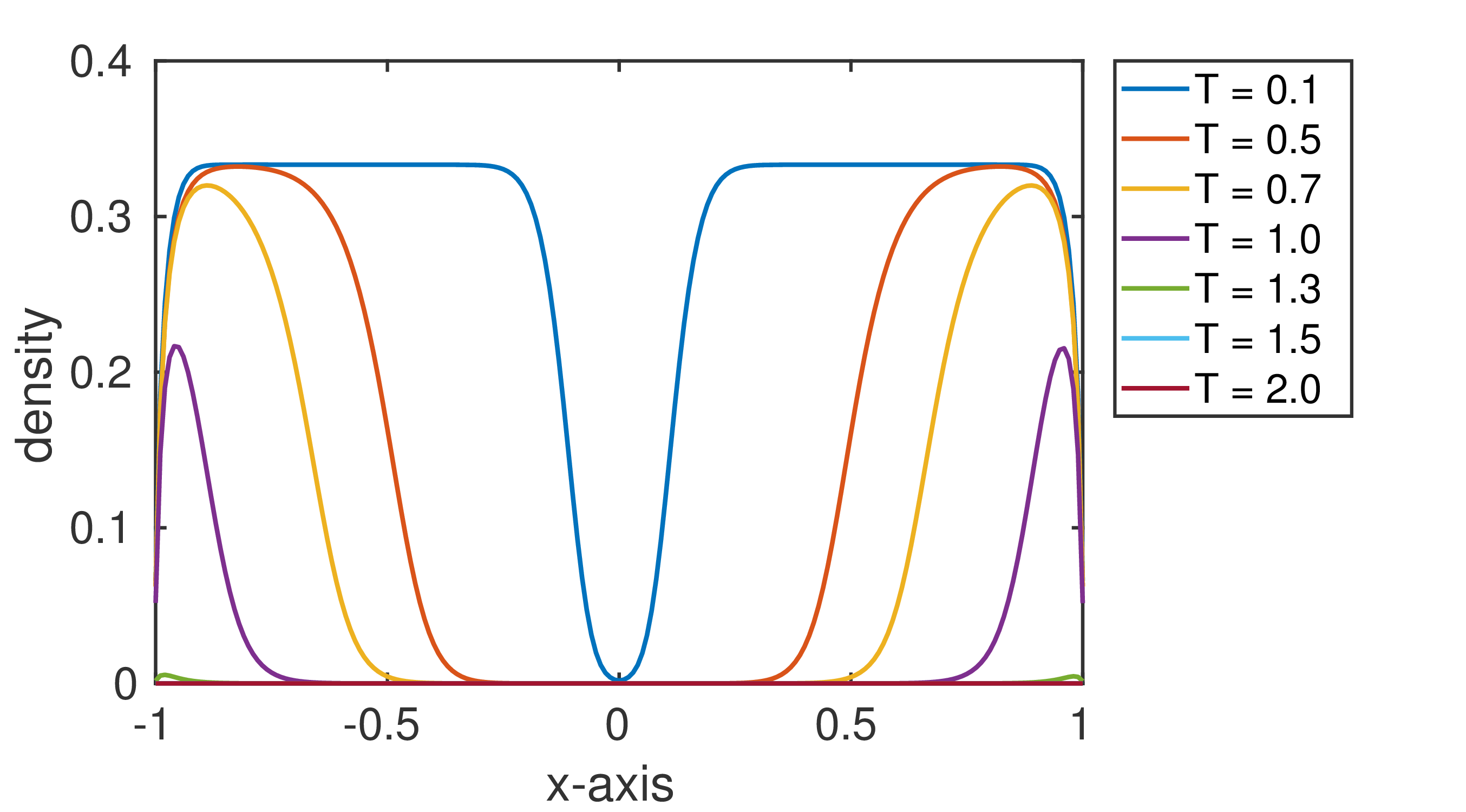}\\ 
	\includegraphics[width=.5\linewidth]{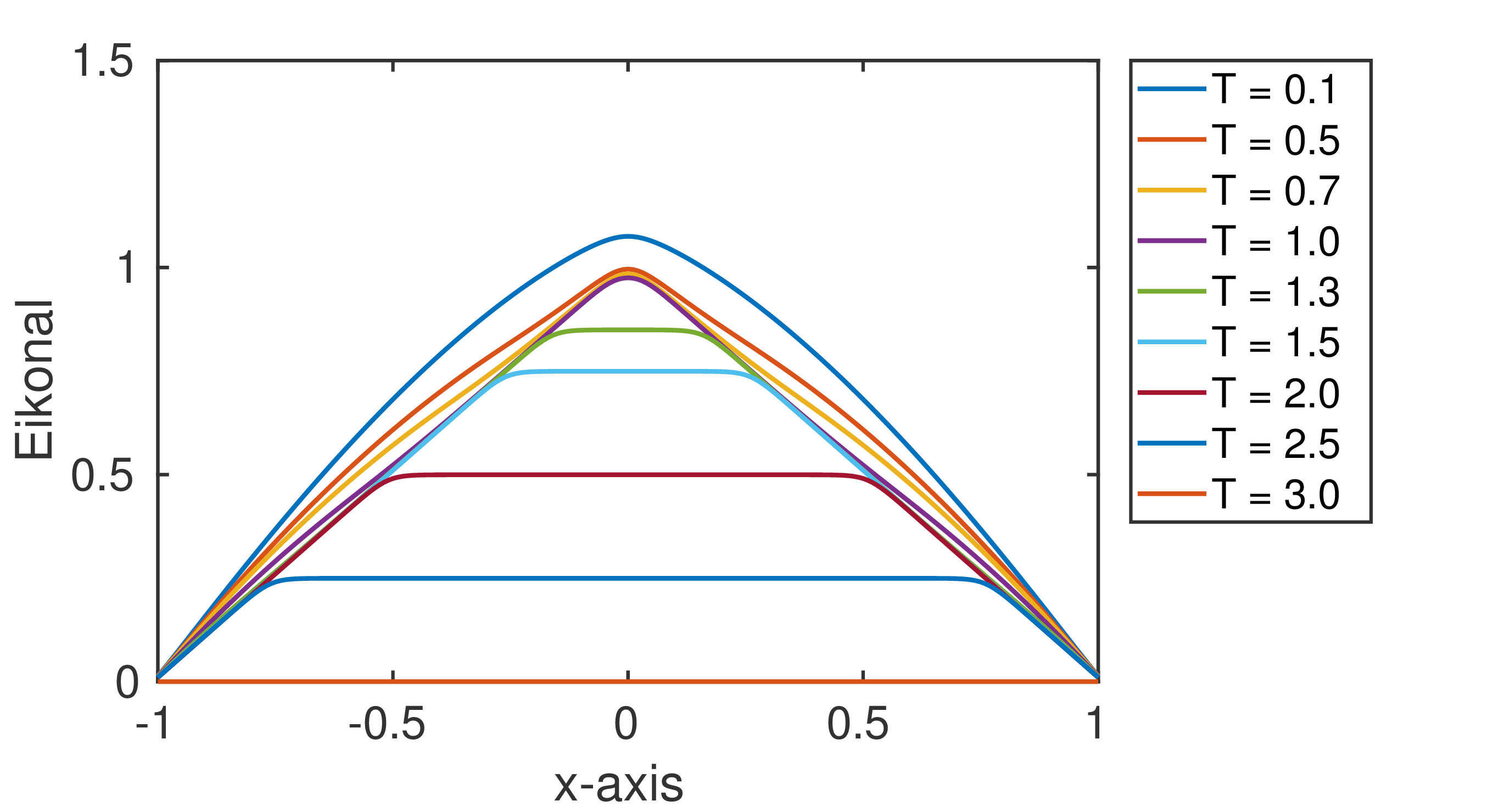}\hfill
	\includegraphics[width=.5\linewidth]{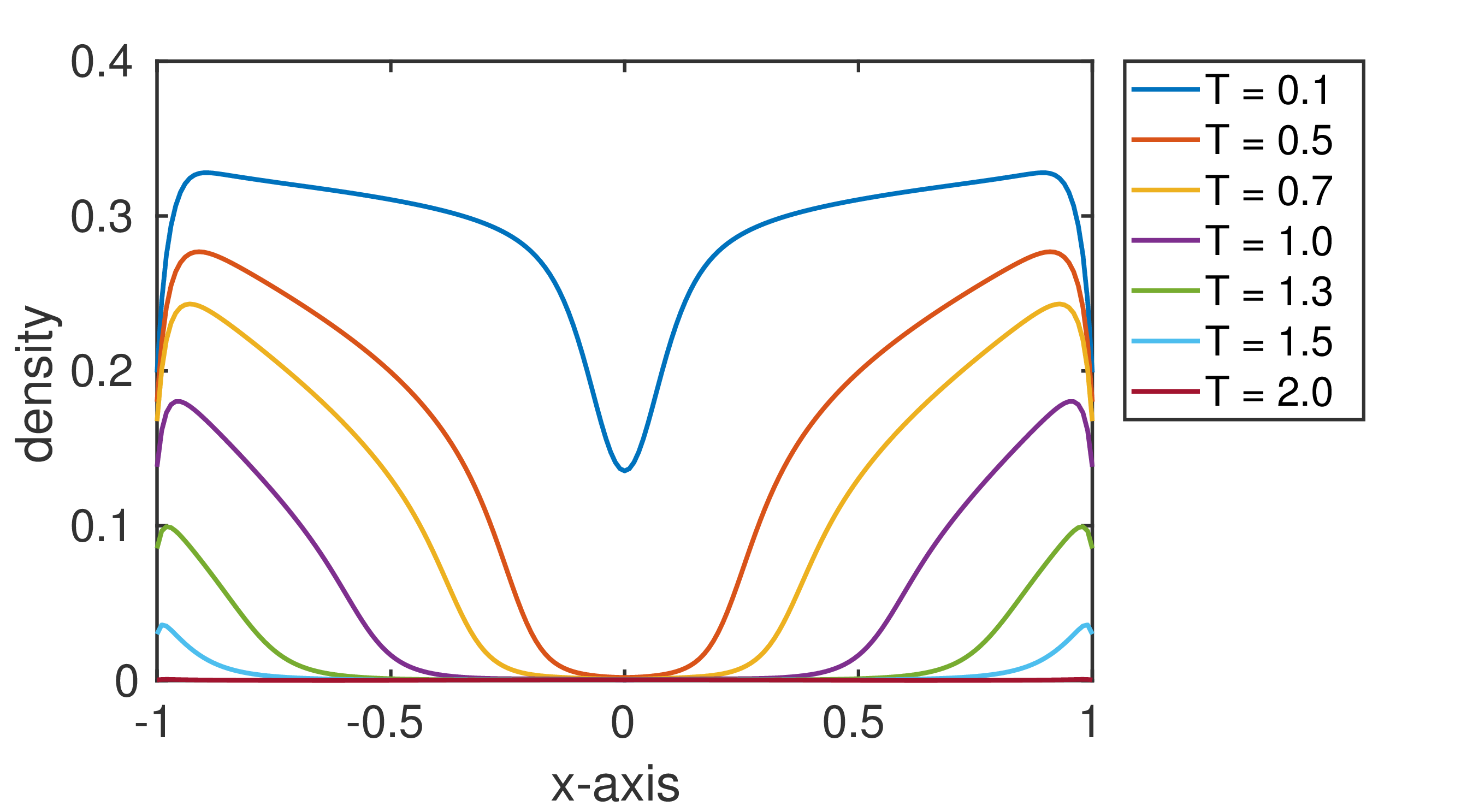}
	\caption{Evolution of solutions through the center along the x-axis at different times for the Hughes model (top) and for the MFG structure (bottom).}\label{fig3}
	
\end{figure}        
Figure \ref{fig3} shows the evolution of solutions through the center along the x-axis for both models at different times. One observes that the solutions in two dimensional case has a similar behavior as the solutions in one dimensional case, see Figure \ref{fig1}.

\section*{Acknowledgment}
	This work is supported by the German research foundation, DFG grant KL 1105/27-1, by RTG GrK 1932 “Stochastic Models for Innovations in the Engineering Sciences”, project area P1 and  by the DAAD PhD program MIC "Mathematics in Industry and Commerce".


\end{document}